\documentclass[11pt]{amsart}
\usepackage{amsfonts, amsmath, enumerate, amssymb}
\usepackage[margin=1in]{geometry}

\newtheorem{theorem}{Theorem}[section]

\newtheorem{lemma}[theorem]{Lemma}
\newtheorem{corollary}[theorem]{Corollary}
\newtheorem{proposition}[theorem]{Proposition}
\theoremstyle{definition}
\newtheorem{definition}[theorem]{Definition}
\theoremstyle{notation}

\theoremstyle{remark}
\newtheorem{remark}[theorem]{Remark}

\numberwithin{equation}{section}
\usepackage{amscd}
\usepackage{xypic}
\usepackage{amsmath, amssymb}

\numberwithin{equation}{subsection}
\newcommand{\be}%
  {\protect\setcounter{equation}{\value{subsubsection}}}
  \newcommand{\ee}%
   {\protect\setcounter{subsubsection}{\value{equation}}}

  {\protect\setcounter{subsubsection}{\value{equation}}}

\def \rmA{\rm A}
\def \rmDA{\rm DA}

\def \rmB{rmB}
\def \rmDB{\rm DB}

\def \EG{\rm EG}
\def \BN(T){\rm BN(T)}

\def \rmB{\rm B}

\def \rmC{\rm C}

\def \Cl{\mathbb C}

\def \colimK.{\underset {\underset K^.  \rightarrow}  {\hbox {lim}}}

\def \colimU.{\underset {\underset U_.  \rightarrow}  {\hbox {lim}}}

\def \rmDX{\rm DX}
\def \rmDY{\rm DY}
\def \rmDZ{\rm DZ}
\def \DS1X{\rm {DS^1X}}
\def \rmD{\rm D}

\def \rmE{{\rm E}}
\def \cE{\mathcal E}

\def \EG1{E{(G \times {\mathbb C}^*)}{\underset {G\times {\mathbb C}^*} 
\times}}

\def \EZ(s)1{E{(Z(s) \times {\mathbb C}^*)}{\underset {(Z(s)\times {\mathbb
C}^*)}  \times}}

\newcommand{\eps}{ \, {\boldsymbol\varepsilon} \,}

\def \EM(u){EM(u){\underset {M(u)}  \times}}
\def \EM(us){EM(u,s){\underset {M(u, s)}  \times}}

\def \EG{\rm EG}

\def \rmf{\rm f}
\def \rmF{\rm F}

\def \g{\it g}

\def \rmG{\rm G}

\def \rmg{\rm g}
\def \GW{\rm GW}

\def \rmh{\rm h}

\def \invlim1{\underset {\infty \leftarrow q}  {\hbox {lim}}^1}

\def \rmj{\rm j}

\def \rmk{\rm k}

\def \k{\it k}

\def \L3{\Lambda \times \Lambda \times \Lambda}
\def \L2{\Lambda \times \Lambda}

\def \lim{\underset \leftarrow  {\hbox {lim}}}

\def \longright2arrow{{\overset \longrightarrow  {\overset {} 
\longrightarrow}}}

\def \L{L\times \Cl ^*}

\def \rmM{\rm M}

\def \cN{\mathcal N}

\def \NT{\rm N_{\rmG}(T)}
\def \rmN{\rm N}

\def \rmp{\rm p}

\def \Spt{\rm {Spt}}

\def \rmP{\rm P}
\def \rmp{\rm p}

\def \rmq{\rm q}

\def \ra{\rightarrow}

\def \RG^{R(G)^{\hat {}}\ }

\def \res{respectively}
\def \rmS{\rm S}
\def \rmR{\rm R}

\def \RHom{{{\mathcal R}{\mathcal H}om}}

\def \rmS{\rm S}

\def\Spt{\rm {\bf Spt}}

\def \SH{{\mathcal S}{\mathcal H}}
\def\Spt{\rm {\bf Spt}}

\def \mbS{\mathbb S}
\def \S1X{\rm {S^1X}}

\def \topGcoh*{^{top, *} _{G}}
\def \topGho*{ _{top,*} ^{G}}

\def \T{{\mathbf T}}

\def \rmT{\rm T}
\def \tr{\it {tr}}
\def \Th{\rm Th}

\def \rmU{\rm U}

\def \rmU{\rm U}

\def \rmX{\rm X}

\def \rmY{\rm Y}

\def \Z(s){Z(s) \times {\mathbb C}^*}
\def \Z{\mathcal Z}

\def \rmZ{\rm Z}

\begin{document}

\title{Additivity of the Motivic Trace and the Motivic Euler-characteristic}

\author{Roy Joshua}
\address{Department of Mathematics, Ohio State University, Columbus, Ohio,
43210, USA.}

\email{joshua.1@math.osu.edu}
\author{Pablo Pelaez}
\address{Instituto de Matem\'aticas, Ciudad Universitaria, UNAM, DF 04510, M\'exico.}
\email{pablo.pelaez@im.unam.mx}
\thanks{\indent Both authors would like to thank the Isaac Newton Institute for Mathematical Sciences, Cambridge, for support and hospitality during the programme {\it K-Theory, Algebraic Cycles and Motivic Homotopy Theory} where part of the work on this paper was undertaken. This work was supported by EPSRC grant no EP/R014604/1.\\ 
	}
\begin{abstract} In this paper,  we settle an open conjecture  regarding the assertion that
	the Euler-characteristic of $\rmG/\NT$ for a split reductive group scheme $\rmG$ and the normalizer of a split maximal torus $\NT$ over a field is $1$ in 
	the Grothendieck-Witt ring with the characteristic exponent of the field inverted, under the assumption that the base field contains a $\sqrt -1$. Numerous applications of this to splittings in the motivic stable homotopy category and to Algebraic K-Theory are worked out in several related papers by Gunnar Carlsson and the
	authors.
\vskip .1cm
	
\end{abstract}
\maketitle

\centerline{\bf Table of contents}
\vskip .2cm 
1. Introduction
\vskip .2cm
2. Additivity of the motivic trace
\vskip .2cm
3. Proofs of the main theorems
\vskip .2cm
\markboth{ Roy Joshua and Pablo Pelaez}{Additivity of the Motivic Trace and the Motivic Euler-characteristic}
\input xypic

\vfill \eject
\section{\bf Introduction}
\label{intro}
This paper is a continuation of the earlier work, \cite{CJ20}, where Carlsson and one of the present authors 
set up motivic and \'etale variants of the classical Becker-Gottlieb transfer. If one recalls, the power and utility of the classical Becker-Gottlieb transfer 
stems from the fact it provided a convenient mechanism to obtain splittings to certain  maps
in the stable homotopy category.  In view of the fact that the transfer is a map of spectra, it induces a map of
the Atiyah-Hirzebruch spectral sequences associated to generalized cohomology theories, and reduces the question on
the existence of stable splittings to the calculation of certain Euler-characteristics in singular cohomology.
The most notable example of this is the calculation of the Euler-characteristic of $\rmG/\NT$, where $\rmG$ is a compact Lie group
and $\NT$ is the normalizer of a maximal torus in $\rmG$. Using the transfer, it becomes then possible to show that
the generalized cohomology of the Borel construction with respect to $\rmG$ for a space $\rmX$ acted on by $\rmG$, is a
split summand of the generalized cohomology of the corresponding Borel construction for $\rmX$ with respect to $\NT$.
This then provided numerous applications, such as double coset formulae for actions of compact groups, generalizing the well-known
double coset formulae for the action of finite groups: see \cite{BG75}, \cite{BG76}, \cite{Beck74}.
\vskip .1cm
The motivic analogue of the statement that the Euler characteristic of $\rmG/\NT$ is $1$
in singular cohomology for compact Lie groups is a conjecture due to Morel (see the next page for more details), that a suitable motivic Euler characteristic
in the Grothendieck-Witt group is $1$, for $\rmG/\NT$, where $\rmG$ is a split connected reductive group and $\NT$ is the 
normalizer of a maximal torus in $\rmG$. We provide an affirmative solution to this conjecture in this paper assuming that the base field $k$ contains a $\sqrt -1$, the precise details of
which are discussed below.
\vskip .2cm
 Let  $k$ denote a {\it perfect field of arbitrary characteristic: we will restrict to the category of smooth quasi-projective schemes over $k$ and adopt the framework of \cite{MV}.}
 Throughout, $\T$ will denote ${\mathbb P}^1$ pointed by $\infty$ and $\T^n$ will denote 
$\T^{\wedge n}$ for any integer $n \ge 0$. $\mbS_{\k}$ will denote the corresponding motivic sphere spectrum. Let $\Spt(\k_{mot})$  denote
 the  category of motivic spectra over $\k$.
 The corresponding stable homotopy category will be denoted  ${\mathcal S}{\mathcal H}(k)$. In positive characteristic ${\rm p}$,  we 
 consider $\Spt(\k_{mot})[{\rm p}^{\rm -1}]$: we will identify this with the category of motivic spectra that are module spectra over the localized 
 sphere spectrum $\mbS_{\k}[{\rm p}^{\rm -1}]$. Then assuming $char(\k)= {\rm 0}$, given a smooth scheme
$\rmX$ of finite type over $k$, $\Sigma_{\T}^{\infty}\rmX_+$ denotes the $\T$-suspension spectrum of $\rmX$ and 
$\rmD(\Sigma_{\T}^{\infty}\rmX_+) = \RHom(\Sigma_{\T}^{\infty}\rmX_+, \mbS_{\k})$, where $\RHom$ denotes the derived internal hom in the category $\Spt(\k_{mot})$.
When $char(\k)= {\rm p} > {\rm 0}$, $\rmD(\Sigma_{\T}^{\infty}\rmX_+) = \RHom(\Sigma_{\T}^{\infty}\rmX_+, \mbS_{\k}[{\rm p}^{-1}])$.
 $\rmD(\Sigma_{\T}^{\infty}\rmX_+)$ is the
{\it Spanier-Whitehead dual} of $\Sigma_{\T}^{\infty}\rmX_+$. It is known (see \cite{Ri13} and also \cite{Ri05}) that after inverting the characteristic exponent, $\Sigma_{\T}^{\infty}\rmX_+$ is
{\it dualizable} in the sense that the natural map $\Sigma_{\T}^{\infty}\rmX_+ \ra \rmD ( \rmD(\Sigma_{\T}^{\infty}\rmX_+))$ is an isomorphism in ${\mathcal S}{\mathcal H}(k)$.
\vskip .1cm
In this context, we have the {\it co-evaluation map} 
\[\mbS_{\k} {\overset c \ra} \Sigma_{\T}^{\infty}\rmX_+ \wedge \rmD(\Sigma_{\T}^{\infty}\rmX_+)\]
and the 
{\it evaluation map} 
\[\rmD(\Sigma_{\T}^{\infty}\rmX_+) \wedge \Sigma_{\T}^{\infty}\rmX_+ \ra \mbS_{\k}\]
in characteristic $0$. In positive characteristic $p$, we also have the co-evaluation map  
\[\mbS_{\k}[\rmp^{-1}] {\overset c \ra} \Sigma_{\T}^{\infty}\rmX_+[\rmp^{-1}] \wedge \rmD(\Sigma_{\T}^{\infty}\rmX_+[\rmp^{-1}])\]
and the 
{\it evaluation map} 
\[\rmD(\Sigma_{\T}^{\infty}\rmX_+[\rmp^{-1}]) \wedge \Sigma_{\T}^{\infty}\rmX_+[\rmp^{-1}] \ra \mbS_{\k}[\rmp^{-1}]. \]
See \cite[p. 87]{DP}. Let $\rmf: \rmX \ra \rmX$ denote a self-map. 
  \vskip .1cm
  \begin{definition} Assume the above setting. Then, in characteristic $0$, the following composition in ${\mathcal S}{\mathcal H}(k)$ 
  defines the {\it trace} ${\tau_{\rmX}(\rmf_+)}$:
   \label{trace}
 \be \begin{equation}
     \label{gen.trace}
    \mbS_{\k} {\overset c \ra} \Sigma_{\T}^{\infty}\rmX_+ \wedge \rmD(\Sigma_{\T}^{\infty} \rmX_+) {\overset {\tau} \ra} \rmD(\Sigma_{\T}^{\infty}\rmX_+) \wedge \Sigma_{\T}^{\infty}\rmX_+  {\overset {id \wedge \rmf} \ra} \rmD(\Sigma_{\T}^{\infty}\rmX_+) \wedge \Sigma_{\T}^{\infty}\rmX_+ {\overset e \ra} \mbS_{\k}.
\end{equation} \ee
\vskip .2cm \noindent
In positive characteristic ${\rm p}$, the following composition in ${\mathcal S}{\mathcal H}(k)[{\rm p}^{-1}]$ defines the corresponding trace,
which will be denoted $\tau_{\rmX, \mbS_{\k}[{\rm p}^{-1}]}(\rmf_+)$:
\be \begin{multline}
     \label{gen.trace}
     \begin{split}
    \mbS_{\k}[\rmp^{-1}] {\overset c \ra} \Sigma_{\T}^{\infty}[\rmp^{-1}]\rmX_+ \wedge \rmD(\Sigma_{\T}^{\infty}[\rmp^{-1}] \rmX_+) {\overset {\tau} \ra} \rmD(\Sigma_{\T}^{\infty}[\rmp^{-1}]\rmX_+) \wedge \Sigma_{\T}^{\infty}[\rmp^{-1}]\rmX_+ \\
    {\overset {id \wedge \rmf} \ra} \rmD(\Sigma_{\T}^{\infty}[\rmp^{-1}]\rmX_+) \wedge \Sigma_{\T}^{\infty}[\rmp^{-1}]\rmX_+ {\overset e \ra} \mbS_{\k}[\rmp^{-1}].
    \end{split}
\end{multline} \ee
\vskip .2cm \noindent
Here $\tau$ is the map interchanging the two factors.  When $\rmf=id_{\rmX}$, the corresponding trace $\tau_{\rmX_+} = \tau_{\rmX}(id_{X+})$ ($\tau_{\rmX_+, \mbS_{\k}[\rmp^{-1}]} = \tau_{\rmX, \mbS_{\k}[\rmp^{-1}]}(id_{\rmX_+})$)  will be denoted $\chi_{mot}(\rmX)$ and
{\it called the motivic Euler-characteristic} of $\rmX$.
  \end{definition}
By \cite{Mo4} (see also \cite{Mo12}), $\pi_{0, 0}(\mbS_{\k})$ identifies with the Grothendieck-Witt ring of the field $k$, $\GW({\it k})$, and therefore  $\chi_{mot}(\rmX)$
 is a class in $\GW({\it k})$ in characteristic $0$ and in $\GW({\it k})[\rmp^{-1}]$ in positive characteristic. 
 (In \cite{Mo4} the isomorphism of $\pi_{0, 0}(\mbS_{\k})$ with the Grothendieck-Witt ring of the field $k$ was proven
 under the assumption $char(\k) \ne 2$: the above restriction on the characteristic of the field ${\it k}$ is removed in 
	 \cite[Theorem 10.12]{BH}.)
 \vskip .2cm
 Then an open conjecture in this setting (due to Morel (see \cite{Lev18})) was the following: let $\rmG$ denote a split reductive group over $k$, with $\rmT$ a split maximal torus and
 $\NT$ its normalizer in $\rmG$. Then the conjecture states that $\chi_{mot}(\rmG/\NT) =1$ in $\GW({\it k})$ with the characteristic exponent of the field inverted. In fact, this is the {\it strong form}
 of the conjecture. The {\it weak form} of the conjecture is simply the statement that $\chi_{mot}(\rmG/\NT)$ is a {\it unit} in $\GW({\it k})$ with the  characteristic exponent of the field inverted.
 \vskip .2cm
 The main result of the current paper is an affirmative solution of the above conjecture as stated in the following theorem.
 \begin{theorem}
  \label{main.thm}
  Let $\rmG$ denote a {\it split} linear algebraic group over the perfect field $k$, with $\rmT$ a split maximal torus and ${\rm N(T)}$ its normalizer in $\rmG$.
  Then the following are true:
  \begin{enumerate}[\rm(i)]
   \item $\chi_{mot}(\rmG/\NT) =1$ in $\GW({\it k})$ if $char(k)=0$ and $k$ contains a $\sqrt -1$.
   \item $\chi_{mot}(\rmG/\NT) =1$ in $\GW({\it k})[\rmp^{-1}]$ if $char(k)=\rmp>0$ and $k$ contains a $\sqrt -1$.
   \end{enumerate}
\end{theorem}
The statements (i) and (ii)  already were used in the preprint \cite{JP20} in the context of proving the additivity also for the transfer and proving various applications
of these in the motivic stable homotopy category. Then, Ananyevskiy proved independently the {\it weak form} of the conjecture in \cite{An}. 
In this paper, we also show how
to simplify the proof discussed in \cite{An}, by making use of our proof of the strong form of the conjecture for fields that contain $\sqrt -1$. This is discussed in 
our proof of the following Corollary.
\begin{corollary}
  \label{cor}
  Assume as in Theorem ~\ref{main.thm} that  $\rmG$ denotes a {\it split} linear algebraic group over $k$, with $\rmT$ a split maximal torus and ${\rm N(T)}$ its normalizer in $\rmG$.
  Then the following hold:
  \begin{enumerate} [\rm(i)]
  \item $\chi_{mot}(\rmG/\NT)$ in $\GW({\it k})$ is a {\it unit} if $char(k)=0$, and 
   \item $\chi_{mot}(\rmG/\NT)$ in $\GW({\it k})[\rmp^{-1}]$ is a {\it unit} if $char(k)=p>0$.
  \end{enumerate}
\end{corollary}

\vskip .2cm
 
 In view of the interest in these results, and because a proof of additivity for the trace is relatively straight-forward, \footnote{i.e., unlike
  additivity for the transfer which is much more involved, and needs the notion of rigidity in an essential manner 
  } 
  we have decided to write this short paper entirely devoted to a self-contained proof of these results. On feeding this
 result into the motivic variant of the transfer constructed in \cite{CJ20} and \cite{CJP22}, we obtain a number of splitting results. The following
 result should serve as a proto-typical example of such applications. 
 \vskip .2cm
 Let $ \rmE \ra \rmB$ denote a $\rmG$-torsor for the action of {\it any} linear algebraic group $\rmG$ with both $\rmE$ and $\rmB$ smooth quasi-projective schemes over $k$, with $\rmB$ {\it connected}.
 Let $\rmY$ denote a $\rmG$-scheme or an unpointed simplicial presheaf provided with a $\rmG$-action. Let  
 	 $\rmq: \rmE{\underset {\rmG} \times}({\rmG}{\underset {\NT} \times}\rmY ) \ra \rmE{\underset {\rmG} \times} \rmY$ denote the map induced by the map ${\rmG}{\underset {\NT} \times}\rmY \ra \rmY$
 	  sending $(g, y) \mapsto gy$. (In case the group $\rmG$ is {\it special} in the sense of Grothendieck, the quotient construction above can be carried out on the 
 	  Nisnevich site, but in general this needs to be carried out on the etale site and then pushed forward to the Nisnevich site, by  means of a derived push-forward: the details
 	  are in \cite[3.4.2]{CJP22}.)
 	  
 	  Then the induced map
 	 \[\rmq^*: \rmh^{\bullet, *}(\rmE{\underset {\rmG} \times} \rmY, M) \ra   \rmh^{\bullet, *}(\rmE{\underset {\rmG} \times}({\rmG}{\underset {\NT} \times}\rmY ), M)\]
 	 is a split injection, where $\rmh^{*, \bullet}$ denotes any generalized motivic cohomology theory with respect to 
 	the motivic spectrum $\rmM$.
 \vskip .2cm	
In order to show that the map $q^*$ is a split monomorphism, one needs the  {\it the transfer }:
\[\tr({\rmY}):{\rm \Sigma}_{\T}^{\infty}(\rmE\times _{\rm G}\rmY)_+ \ra \Sigma_{\T}^{\infty} (\rmE{\underset {\rmG} \times}({\rmG}{\underset {\NT} \times}\rmY)_+,\]
which is a map in ${\mathcal S}{\mathcal H}(k)$ (${\mathcal S}{\mathcal H}(k)[\rmp^{-1}]$, \res)  so that the composition $\tr({\rmY})^* \circ q^*$
is multiplication by  $\chi_{mot}(\rmG/{\rm N(T)})$. Therefore, knowing $\chi_{mot}(\rmG/{\rm N(T)})$ is a unit in $\GW(\k)$ shows $q^*$ is a split injection.
\vskip .2cm
See also \cite{JP22}, where such splittings obtained from the 
 motivic transfer proves a variant of the classical Segal-Becker theorem for Algebraic K-Theory.
 
\begin{remark} Certain special cases of the above splitting results, for groups that are {\it special} seem to be also worked out in \cite{Lev18}, {\it under
 the assumption the above conjecture is true}. Observe that a linear algebraic group $\rmG$ is {\it special} in the sense of Grothendieck (see \cite{Ch}), if
 any torsor for $\rmG$ is locally trivial on the Zariski site. {\it Special} groups include $\{{\rm GL}_n, {\rm SL}_n|n\}$, but {\it exclude
 all orthogonal groups as well as finite groups}. For groups $\rmG$ that are {\it not special}, $\rmG$-torsors are locally trivial only on the \'etale site.
 The construction of the transfer, worked out in \cite{CJ20} and \cite[Chapter 3]{CJP22} apply for all such groups. 
\end{remark}

 \vskip .2cm
 Here is an overview of the paper. One of the key techniques that is used in the proof of Theorem ~\ref{main.thm} is to show that the trace and the motivic Euler-characteristic are 
 {\it additive up to multiplication by a sign} in general, and additive when the base field $\k$ has a $\sqrt -1$.
 We devote section 2 of the paper to establishing this additivity. Section 3 then completes the proof of the above theorem closely following
  the ideas for a proof of the corresponding result as in \cite[Lemma 3.5]{BP} in the \'etale setting. 
 \vskip .2cm
  {\bf Acknowledgments}. The first author would like to thank Gunnar Carlsson
for getting him interested in the problem of constructing a Becker-Gottlieb transfer in the motivic framework and for numerous helpful discussions. Both authors  would like
to thank Michel Brion for helpful discussions on fixed point schemes as well as on aspects of Theorem ~\ref{torus.act}. We are also
happy to acknowledge \cite[Lemma 3.5 and its proof]{BP} as one of the inspirations for this paper. We also thank Alexey Ananyevskiy 
for helpful comments on our preprint \cite{JP20},  which have enabled us to sharpen our results, and also for bringing his results to our attention. Finally, it is a pleasure to
acknowledge our intellectual debt to Fabien Morel and Vladimir Voevodsky for their foundational work in motivic homotopy theory. In addition, the authors are also grateful to the 
referee for providing us very valuable feedback, which surely have helped us sharpen some results and improve the overall organization.

 \section{\bf Additivity of the Motivic Trace}
The main goal of this section is to establish additivity properties for the pre-transfer and trace. 
But we begin by establishing certain properties of a general nature for the pre-transfer and the trace. 
\subsection{Basic properties of the pre-transfer and trace}
 It is convenient to reformulate the trace in terms of the pre-transfer, which we proceed to discuss next. At the same time, we extend the framework as follows. The following 
discussion is a variant of what appears in \cite[Chapter III]{LMS}. See also \cite{May01}, \cite{GPS} and \cite{HSS} for related discussions.
 \begin{definition} 
\label{co-mod.st}
(Co-module structures)
Assume that $\rmC$ is an unpointed simplicial presheaf, i.e., $\rmC$ is a contravariant functor from a given site to 
the category of unpointed simplicial sets. Let $\rmC_+$ denote the corresponding pointed simplicial presheaf. Then the diagonal map $\Delta : \rmC_+ \ra \rmC_+ \wedge \rmC_+$ together with
the augmentation $\epsilon: \rmC_+ \ra \rmS^0$ defines the structure of an associative co-algebra of simplicial presheaves on $\rmC_+$.
A pointed simplicial presheaf $\rmP$ will be called a right $\rmC_+$-co-module, if it comes equipped with maps $\Delta: \rmP \ra \rmP \wedge \rmC_+$
so that  the diagrams:
\begin{equation}
 \label{diagonal.1}
  \xymatrix{{ \rmP} \ar@<1ex>[r]^{\Delta} \ar@<1ex>[d]^{\Delta} & {\rmP \wedge \rmC_+} \ar@<1ex>[d]^{id \wedge \Delta}\\  
             {\rmP \wedge \rmC_+} \ar@<1ex>[r]^{\Delta \wedge id} & {\rmP \wedge \rmC_+ \wedge \rmC_+ }}
  \quad \mbox{ and } \xymatrix{ {\rmP} \ar@<1ex>[d]^{\Delta} \ar@<1ex>[dr]^{id} \\           
                                {\rmP \wedge \rmC_+} \ar@<1ex>[r]^{id \wedge \epsilon} & {\rmP \wedge \rmS^0}}
\end{equation}
commute. {\it The most common choice of $\rmP$ is with $\rmP = \rmC_+$} and with the obvious diagonal map $\Delta: \rmC_+ \ra \rmC_+ \wedge \rmC_+$ as 
providing the co-module structure. However, the reason we are constructing the pre-transfer in this generality (see the definition below) is so that we are able to obtain
strong additivity results as in Theorem ~\ref{additivity.tr.0}. 
\end{definition}
 
 \vskip .1cm
 \begin{definition} ({\it The  pre-transfer})
 	\label{pretransfer}
 	Assume that the pointed simplicial presheaf $\rmP$ is such that $\Sigma_{\T}^{\infty} \rmP$ is dualizable in $\Spt(\k_{mot})$ and is provided
 	with a map $\rmf: \rmP \ra \rmP$. Assume further that $\rmC$ is an unpointed simplicial presheaf so that $\rmP$ is a
 	right $\rmC_+$-co-module. 
 	Then the {\it  pre-transfer with respect to $\rmC_+$} is defined to be a map
 	$tr'(\rmf): \mbS_{\k} \ra \Sigma_{\T}^{\infty} \rmC_+$, which is the
 	composition of the following maps. Let $e: \rmD(\Sigma_{\T}^{\infty} \rmP )
 	\wedge \Sigma_{\T}^{\infty}  \rmP \ra \mbS_{\k}$ denote the evaluation map. 
 	We take the dual of this map to obtain:
 	\be \begin{equation}
 	  \label{coeval}
 	c=\rmD(e): \mbS_{\k} \simeq \rmD(\mbS_{\k}) \ra \rmD(\rmD(\Sigma_{\T}^{\infty} \rmP ) \wedge (\Sigma_{\T}^{\infty}  \rmP) ) {\overset {\simeq} \leftarrow} \rmD(\Sigma_{\T}^{\infty}  \rmP)\wedge 
 	(\Sigma_{\T}^{\infty} \rmP ) {\overset {\tau} \ra}  (\Sigma_{\T}^{\infty} \rmP ) \wedge \rmD(\Sigma_{\T}^{\infty} \rmP).
 	\end{equation}
 	
 	Here $\tau$ denotes the obvious flip map interchanging the two factors and $c$ denotes the co-evaluation. The reason that taking
 	the double dual yields the same object up to weak-equivalence
 	is because we are in fact taking the dual in the setting discussed
 	above.  Observe that all the maps that {\it go in the left-direction are weak-equivalences}. All the maps involved in the definition of the co-evaluation map are  {\it natural maps}.

 	To complete the definition of the pre-transfer, one simply composes
 	the co-evaluation map with the following composite map:
 	\be \begin{multline}
 	\begin{split}
 	(\Sigma_{\T}^{\infty}  \rmP )  \wedge \rmD(\Sigma_{\T}^{\infty} \rmP) {\overset {\tau} \ra} \rmD(\Sigma_{\T}^{\infty} \rmP)\wedge (\Sigma_{\T}^{\infty} \rmP ){\overset {id \wedge f} \ra } \rmD(\Sigma_{\T}^{\infty} \rmP)\wedge (\Sigma_{\T}^{\infty} \rmP ) \\
 	 {\overset {id \wedge \Delta} \ra }\rmD(\Sigma_{\T}^{\infty}  \rmP)\wedge (\Sigma_{\T}^{\infty}  \rmP ) \wedge (\Sigma_{\T}^{\infty} \rmC_+ )
 	{\overset {e \wedge id} \ra } \mbS_{\k} \wedge (\Sigma_{\T}^{\infty}  \rmC_+ ) \simeq \Sigma_{\T}^{\infty}  \rmC_+.
 	\end{split}
 	\end{multline} \ee
 	\vskip .1cm
 	The corresponding {\it trace} $\tau(\rmf)$, is defined as the composition of the above pre-transfer $tr'(\rmf)$ with the projection $\pi$ sending 
 	$\rmC_+$ to $\rmS^0_+$. 
 	\vskip .1cm
 	When  $\rmf=id_{\rmP}$, the pre-transfer (trace) will be denoted $tr_{\rmP}'$ ($\tau_{\rmP}$, \res), and when
 	$\rmP =\rmC_+ $ and $\rmf=id_{\rmP}$, the pre-transfer (trace) will be denoted $tr_{\rmC_+}'$ ($\tau_{\rmC_+}$, \res). 
 \end{definition}
 \begin{remark}
  \label{key.ident.trace} Observe that now the trace map $\tau_{\rmC_+}$ identifies with the following composite map:
   \[  \tau_{\rmC_+}: \mbS_{\k} {\overset c \ra} \Sigma_{\T}^{\infty}  \rmC_+ \wedge \rmD(\Sigma_{\T}^{\infty} \rmC_+) {\overset {\tau} \ra} \rmD(\Sigma_{\T}^{\infty}  \rmC_+) \wedge \Sigma_{\T}^{\infty}  \rmC_+ {\overset e \ra} \mbS_{\k} .\]
 \end{remark}
 \begin{definition}
 \label{generalized.traces}
 If $\cE$ denotes {\it any commutative ring spectrum} in $\Spt(\k_{mot})$, for example, $\mbS_{\k}[\rmp^{-1}]$, we will let
 $\Spt(\k_{mot}, \cE)$ denote the category of $\cE$-module spectra over $\cE$.
 Then one may replace the sphere spectrum
$\mbS_{\k}$ everywhere by $\cE$ and define the pre-transfer and trace similarly, provided the
 unpointed simplicial presheaf $\rmC$ is such that $\cE \wedge \rmC_+$ is dualizable in ${\Spt}(\k_{mot}, \cE)$ and is provided
 	with a  map $\rmf: \rmC \ra \rmC$.  When $\rmP = \rmC_+$, these will be
 denoted $\tr(\rmf_+)_{\cE}'$,  $\tr'_{\rmC_+, \cE}$, $\tau_{\rmC_+, \cE}$, etc.
 \end{definition}

 \vskip .2cm
 Let 
 \be \begin{equation}
      \label{additivity.0}
 \rmU_+ {\overset {\rmj_+} \ra} \rmX_+ {\overset {\rmk_+} \ra} \rmX/\rmU = Cone(\rmj_+) \ra \rmS^1 \wedge \rmU_+
 \end{equation} \ee
 denote a cofiber sequence where both $\rmU$ and $\rmX$ are unpointed simplicial presheaves, with $j_+$ a cofibration. Now a key point to observe is that all of $\rmU$, $\rmX$ and $\rmX/\rmU$ have the
 structure of right $\rmX_+$-co-modules. The right $\rmX_+$-co-module structure on $\rmX_+$ is given by the diagonal map $\Delta: \rmX_+ \ra \rmX_+ \wedge \rmX_+$,
 while the right $\rmX_+$-co-module structure on $\rmU_+$ is given by the map $\Delta: \rmU _+ {\overset {\Delta} \ra} \rmU_+ \wedge \rmU_+ {\overset {id \wedge j_+} \ra} \rmU_+ \wedge \rmX_+$,
 where $j: \rmU \ra \rmX$ is the given map. The right $\rmX_+$-co-module structure on $\rmX/\rmU$ is obtained in view of the commutative square
 \vskip .1cm
 \be \begin{equation}
      \xymatrix{{\rmU} \ar@<1ex>[r]^{(id \times j)\circ \Delta} \ar@<1ex>[d]^j & {\rmU \times X} \ar@<1ex>[d]^{j \times id}\\
                {\rmX} \ar@<1ex>[r]^{\Delta} & {\rmX \times \rmX}}
 \end{equation} \ee
 which provides the map 
 \be \begin{equation}
   \label{comod.XU}
 \rmX/\rmU \ra (\rmX \times \rmX)/(\rmU \times \rmX) \cong (\rmX/\rmU) \wedge \rmX_+.
\end{equation} \ee

 \vskip .1cm 
We begin with the following results, which are variants of \cite[Theorem 7.10, Chapter III and Theorem 2.9, Chapter IV]{LMS} adapted to our contexts.
 \begin{theorem} 
  \label{additivity.tr.0}
  Let $\rmU_+ {\overset {\rmj_+} \ra} \rmX_+ {\overset {\rmk_+} \ra} \rmX/\rmU = Cone(\rmj) \ra \rmS^1 \wedge \rmU_+$ denote a cofiber sequence 
  as in ~\eqref{additivity.0}. Let $f: \rmU_+ \ra \rm U_+$, $g: \rmX_+ \ra \rmX_+$ denote two pointed maps so that the diagram
  \vskip .1cm
  \[\xymatrix{{\rmU_+} \ar@<1ex>[r]^{j_+} \ar@<1ex>[d]_{f} & {\rmX_+} \ar@<1ex>[d]_{g} \\
            {\rmU_+} \ar@<1ex>[r]^{j_+} & {\rmX_+} } \]
 \vskip .1cm \noindent
 commutes. Let  $h: \rmX/\rmU \ra \rmX/\rmU$ denote the corresponding induced map. Then, with the right $\rmX_+$-co-module structures discussed above,
 one obtains the following commutative diagram:
 \be \begin{equation}
      \label{additivity.tr.diagram.0}
\xymatrix{{\rmU_+} \ar@<1ex>[r]^{j_+} \ar@<1ex>[d]_{\Delta} & {\rmX_+} \ar@<1ex>[r]^{k_+} \ar@<1ex>[d]_{\Delta} & {\rmX/\rmU} \ar@<1ex>[r]^l \ar@<1ex>[d]_{\Delta} & {\rmS^1\wedge \rmU_+} \ar@<1ex>[d]_{\rmS^1 \Delta}\\
            {\rmU_+\wedge \rmX_+} \ar@<1ex>[r]^{j_+ \wedge id} & {\rmX_+ \wedge \rmX_+} \ar@<1ex>[r]^{k_+ \wedge id} & {(\rmX/\rmU) \wedge \rmX_+} \ar@<1ex>[r]^{l \wedge id}  & {\rmS^1\wedge \rmU_+ \wedge \rmX_+}.}
 \end{equation} \ee
 \vskip .1cm \noindent
  Assume further that the ${\T}$-suspension spectra of all the above simplicial presheaves are dualizable in 
  $\Spt(\k_{mot})$. Then, one obtains in $\SH(k)$:
  \[\tr'(g) = \tr'(f) + \tr'(h),\quad \mbox{  and } \tau(g) = \tau(f) + \tau(h).\]
 
 \vskip .1cm
 Let $\cE$ denote a commutative ring spectrum in ${\Spt}(\k_{mot})$. 
  Then the 
 	corresponding results also hold if the smash products of the above simplicial presheaves
 	with the ring spectrum $\cE$ are dualizable in ${\Spt}(\k_{mot}, \cE )$.
 
 \end{theorem}

 \begin{theorem}
 	\label{additivity.tr.1}
 	Let $\rmF= \rmF_1\sqcup_{\rmF_3}\rmF_2$ denote a  pushout of unpointed 
 	simplicial presheaves on the big Nisnevich site of the base scheme, with 
 	the corresponding maps $\rmF_3 \ra \rmF_2$, $\rmF_3 \ra \rmF_1$ and $\rmF_j \ra \rmF$, for $j=1,2, 3$, assumed to be
 	cofibrations (that is, injective maps of presheaves).
  Assume further the following: the $\T$-suspension spectra of all the above simplicial presheaves are dualizable in $\Spt(\k_{mot})$.
 	Let $i_j:\rmF_j \ra \rmF$ denote the inclusion $\rmF_j \ra \rmF$, $j=1,2, 3$.
 	Then,  one obtains in $\SH(k)$:
 	\begin{enumerate}[\rm(i)]
 	\item 	$tr_{\rmF_+}' =  i_1 \circ tr_{\rmF_{1+}}' +  i_2 \circ tr_{\rmF_{2+}}' -  i_3 \circ tr_{\rmF_{3+}}' \mbox{ and } \tau_{\rmF_+} =  \tau_{\rmF_{1+}} +  \tau_{\rmF_{2+}} -  \tau_{\rmF_{3+}}$, \\
 	where \, $tr_{\rmF_+}'$ and $tr_{\rmF_{j+}}'$, $j=1, 2, 3$ ($\tau_{\rmF_+}$, $\tau_{\rmF_{j+}}$, $j=1,2, 3$) denote the  pre-transfer maps ( trace maps, respectively).
 	\item
 	In particular, taking $\rmF_2=*$, and $\rmF= Cone(\rmF_3 \ra \rmF_1)$, we obtain in $\SH(k)$:\\
	    $tr_{\rmF}' = i_1 \circ tr_{\rmF_{1+}}' - i_3 \circ tr_{\rmF_{3+}}' \mbox{ and } \tau_{\rmF} =  \tau_{\rmF_{1+}} - \tau_{\rmF_{3+}}$.
        \end{enumerate}
    \vskip .1cm
 Let $\cE$ denote a commutative ring spectrum in ${\Spt}(\k_{\rm mot})$. 
  Then the 
 	corresponding results also hold if the smash products of the above simplicial presheaves
 	with the ring spectrum $\cE$ are dualizable in ${\Spt}(\k_{\rm mot}, \cE )$.    
        
 \end{theorem}	
 \vskip .1cm
 {\it Our next goal is to provide proofs of these two theorems}. We will discuss the proofs explicitly only for the case of spectra in 
 $\Spt(\k_{{ mot}})$, as the corresponding results readily extend to spectra in $\Spt(\k_{mot}, \cE)$ for a
  commutative ring spectrum $\cE$ in $\Spt(\k_{mot}, \cE)$. The additivity of the trace follows readily from the additivity of the pre-transfer, as it is obtained by composing with the projection 
 $\Sigma_{\T}^{\infty}  \rmX_+ \ra \mbS_{k}$.
 \vskip .1cm
 Since this is discussed in the topological framework in \cite[Theorem 7.10, Chapter III and Theorem 2.9, Chapter IV]{LMS},
 our proof amounts to verifying carefully and in a detailed manner that the same arguments there carry over to our framework. This is possible, largely because the arguments in the
 proof of \cite[Theorem 7.10, Chapter III and Theorem 2.9, Chapter IV]{LMS} depend only
 on a theory of Spanier-Whitehead duality in a symmetric monoidal triangulated category
  framework and \cite{DP} shows that the entire theory of Spanier-Whitehead duality works in 
  such general frameworks. Nevertheless, it seems prudent to  show explicitly that
  at least the key arguments in \cite[Theorem 7.10, Chapter III and Theorem 2.9, Chapter IV]{LMS} carry over to our framework.
  It may be important to point out that the discussion in \cite[Chapters III and IV]{LMS} is carried out in the equivariant
  framework: as all our discussion is taking place with no group actions, one may take the group to be trivial in the discussion in {\it op. cit}.
 \vskip .2cm
 The very first observation is that the hypotheses of Theorem ~\ref{additivity.tr.0} readily imply the commutativity of the 
 diagram:
 \[\xymatrix{{\rmU_+} \ar@<1ex>[r]^{j_+} \ar@<1ex>[d]_f & {\rmX_+} \ar@<1ex>[r]^{k_+} \ar@<1ex>[d]_g & {\rmX/\rmU} \ar@<1ex>[r]^l \ar@<1ex>[d]_h & {\rmS^1 \wedge \rmU_+} \ar@<1ex>[d]_{\rmS^1 {\it f}}\\
            {\rmU_+} \ar@<1ex>[r]^{j_+} & {\rmX_+} \ar@<1ex>[r]^{k_+}  & {\rmX/\rmU} \ar@<1ex>[r]^l  & {\rmS^1 \wedge \rmU_+}.}\]
 \vskip .1cm \noindent
 Next we proceed to verify the commutativity of the diagram ~\eqref{additivity.tr.diagram.0}. Since the first square clearly commutes, it suffices
 to verify the commutativity of the second square. This follows readily in view of the following commutative square of pairs:
 \[\xymatrix{{(\rmX, \phi)} \ar@<1ex>[r] \ar@<1ex>[d]^{\Delta} & {(\rmX, \rmU)} \ar@<1ex>[d]^{\Delta}\\
            {(\rmX \times \rmX, \phi)} \ar@<1ex>[r] & {(\rmX \times \rmX, \rmU \times \rmX)}.}\]
\vskip .1cm \noindent
Observe, as a consequence that we have verified that the hypotheses of \cite[Theorem 7.10, Chapter III]{LMS} are satisfied by the $\Sigma_{\T}^{\infty}$-suspension spectra
of all the simplicial presheaves  appearing in ~\eqref{additivity.tr.diagram.0}. 
 \vskip .2cm
 The next step is to observe that the $\rmF_i$, $i=1,2, 3$ ($\rmF$) in our Theorem ~\ref{additivity.tr.1}, correspond to the $\rmF_i$ ($\rmF$, \res)  in \cite[Theorem 2.9, Chapter IV]{LMS}. Now observe that
 \be \begin{equation}
 \label{cofiber.seq}
 \rmF_{3+} \ra (\rmF_{1} \sqcup \rmF_2)_+ \ra \rmF_+ \ra \rmS^1 \wedge \rmF_{3,+}
 \end{equation} \ee
 \vskip .1cm \noindent
 is a distinguished triangle. Moreover as $\rmF_1 \sqcup F_2$ has a natural map (which we will call $k$) into $\rmF$, 
  there is a commutative diagram:
   \xymatrix{{(\rmF_1 \sqcup \rmF_2)_+} \ar@<1ex>[r]^{k}  \ar@<1ex>[d]^{\Delta} & 	{\rmF_+}\ar@<1ex>[d]^{\Delta}\\
   {(\rmF_1 \sqcup \rmF_2)_+ \wedge \rmF_+}  \ar@<1ex>[r]^{k \wedge id} & {\rmF_+ \wedge \rmF_+}  .}
  \vskip .1cm \noindent
  Then the distinguished triangle ~\eqref{cofiber.seq} provides the commutative diagram:
  \[\xymatrix{{(\rmF_1 \sqcup \rmF_2)_+} \ar@<1ex>[r]^{k}  \ar@<1ex>[d]^{\Delta} &{\rmF_+}\ar@<1ex>[d]^{\Delta} \ar@<1ex>[r] &{\rmS^1 \wedge \rmF_{3,+}} \ar@<1ex>[d]^{\Delta} \ar@<1ex>[r]&{\rmS^1\wedge (\rmF_1 \sqcup \rmF_2)_+} \ar@<1ex>[d]^{\Delta} \\
  	{(\rmF_1 \sqcup \rmF_2)_+ \wedge \rmF_+}  \ar@<1ex>[r]^{k \wedge id} & {\rmF_+ \wedge \rmF_+}  \ar@<1ex>[r] & {(\rmS^1\wedge \rmF_{3,+} )\wedge \rmF_+} \ar@<1ex>[r] &{(\rmS^1 \wedge (\rmF_1 \sqcup \rmF_2)_+) \wedge \rmF_+}}\]
  \vskip .1cm \noindent
  so that the hypotheses of \cite[Theorem 7.10, Chapter III]{LMS} are satisfied with $\rmX$, $\rmY$ and $\rmZ$ there
  equal to the $\Sigma_{\T}^{\infty}$-suspension spectra of $(\rmF_1 \sqcup \rmF_2)_+$, $\rmF_+$ and $\rmS^1 \wedge \rmF_{3,+}$. These arguments, therefore reduce the proof of Theorem ~\ref{additivity.tr.1} to 
   that of Theorem ~\ref{additivity.tr.0}.
  \vskip .2cm
  Therefore, what we proceed to verify is that, then the proof of \cite[Theorem 7.10, Chapter III]{LMS} carries over to our framework.  This will then complete the proof of
  Theorem ~\ref{additivity.tr.0}.
   A key step of this amounts to verifying that the big commutative diagram given on \cite[p. 166]{LMS} carries over to our framework. One may observe that this big diagram is broken up into various sub-diagrams, labeled (I) through (VII) and that it suffices to verify that each of these sub-diagrams commutes up to homotopy. 
 This will prove that additivity holds for the trace.
  \vskip .2cm
  For this, it seems best to follow the terminology adopted in \cite[Theorem 7.10, Chapter III]{LMS}: therefore
  we will let $\rmU_+$ ($\rmX_+$ and $\rmX/\rmU$) in Theorem ~\ref{additivity.tr.0} be denoted $\rmX$ ($\rmY$
  and $\rmZ$, \res) for the remaining part of the proof of Theorem ~\ref{additivity.tr.0}.
  Let $k: \rmX \ra \rmY$ ($i: \rmY \ra \rmZ$  and $\pi: \rmZ \ra \rmS^1 \wedge \rmX$) denote the corresponding maps
  $j_+:\rmU_+ \ra \rmX_+$ ($k_+:\rmX_+ \ra \rmX/\rmU$, and the map $l:\rmX/\rmU \ra \rmS^1 \wedge \rmU_+$) as in Theorem ~\ref{additivity.tr.0}.
  Then the very first step in this
  direction is to verify that the three squares
  \be \begin{equation}
  \label{three.squares}
       \xymatrix{{\rmDY \wedge \rmX } \ar@<1ex>[r]^{id \wedge k} \ar@<1ex>[d]^{\rmD k \wedge id}& {\rmDY \wedge \rmY} \ar@<1ex>[d]^{e} \\
       	 {\rmDX \wedge \rmX} \ar@<1ex>[r]^e & {\mbS_{\k}}}, \xymatrix{{\rmDZ \wedge \rmY } \ar@<1ex>[r]^{id \wedge i} \ar@<1ex>[d]^{\rmD i \wedge id}& {\rmDZ \wedge \rmZ} \ar@<1ex>[d]^{e} \\
       		{\rmDY \wedge \rmY} \ar@<1ex>[r]^e & {\mbS_{\k}}}, \mbox{ and } \xymatrix{{\rmD (S^1\wedge \rmX)\wedge \rmZ} \ar@<1ex>[r]^{id \wedge \pi} \ar@<1ex>[d]^{\rmD \pi \wedge id}& {\rmD (\rmS^1 \wedge \rmX) \wedge (\rmS^1 \wedge X)} \ar@<1ex>[d]^{e} \\
       		{\rmDZ \wedge \rmZ} \ar@<1ex>[r]^e & {\mbS_{\k}}}
\end{equation} \ee
 \vskip .1cm \noindent
 commute up to homotopy. (The homotopy commutativity of these squares is a formal consequence of Spanier-Whitehead duality: see \cite[pp. 324-325]{Sw} for proofs in the classical setting.) As argued on \cite[page 167, Chapter III]{LMS}, the composite 
 $e \circ (\rmD \pi \wedge i): \rmD(\rmS^1\wedge \rmX) \wedge \rmY \ra \mbS_{\k}$ is equal to $ e\circ((id \wedge \pi) \circ (id \wedge i)$ and is therefore the trivial map. Therefore, if $j$ denotes the 
 inclusion of $ \rmDZ \wedge \rmZ$ in the cofiber of $\rmD \pi \wedge i$, one obtains the
 induced map $\bar e: (\rmDZ \wedge \rmZ)/(\rmD(S^1\wedge \rmX) \wedge \rmY) \ra \mbS_{\k}$ so that the triangle
 \be \begin{equation}
      \xymatrix{{\rmDZ \wedge \rmZ}  \ar@<1ex>[r] ^e \ar@<1ex>[d]^{j} & {\mbS_{\k} }\\
      	{(\rmDZ \wedge \rmZ)/(\rmD(S^1\wedge X) \wedge \rmY) } \ar@<1ex>[ur]^{\bar e}}
 \end{equation}
 \vskip .1cm \noindent
 homotopy commutes. This provides the commutative triangle denoted (I) in \cite[p. 166]{LMS} there and the commutative triangle denoted (II) there commutes by the second and third commutative squares in ~\eqref{three.squares}. The duals of (I) and (II) are the triangles denoted (I*) and (II*) (on \cite[p. 166]{LMS})
 and therefore, they also commute.
\vskip .2cm
Next we briefly consider the homotopy commutativity of the remaining diagram beginning with the
squares labeled (III), (IV) and (V) in \cite[p. 166]{LMS}. Since the maps denoted $\delta$ are weak-equivalences,
 it suffices to show that these squares homotopy commute when the maps denoted $\delta ^{-1}$
 are replaced by the corresponding maps $\delta$ going in the opposite direction.  Such maps
 $\delta$ appearing there  are all special instances of the following natural map:
 $\delta : \rmDB \wedge \rmA \ra D(\rmDA \wedge \rmB)$, for two spectra $\rmA$ and $ \rmB$ in  $\Spt(\k_{mot})$. The homotopy commutativity of the squares (III), (IV) and (V) are reduced therefore to
 the naturality of the above map in the arguments $\rmA$ and $\rmB$: see the discussion in \cite[pp. 167-168]{LMS}. The commutativity of the
 triangle labeled (VI) follows essentially from the definition of the maps there. Finally the homotopy commutativity of the square (VII) is reduced to the following lemma, which is 
 simply a restatement of \cite[Lemma 7.11, Chapter III]{LMS}. These will complete the proof
 for the additivity property for the trace and hence the proofs of Theorems ~\ref{additivity.tr.0} and ~\ref{additivity.tr.1}. 
 \vskip .2cm
\begin{lemma}
   \label{[lemma.LMS]}
   Let $\rmf: \rmA \ra \rmX$ and $\rmg: \rmB  \ra  \rmY$
   be maps in $\Spt(\k_{mot})$ and let
   $i: \rmX \ra  Cone(\rmf)$ and $j : \rmY \ra  Cone(\rmg)$ be the inclusions into their cofibers.
   Then the boundary map $\delta :\Sigma_{\rmS^1}^{-1}Cone(i \wedge j) \ra Cone (\rmf \wedge \rmg)$ in the
   cofiber sequence $Cone(\rmf \wedge \rmg) \ra Cone ((i \circ \rmf) \wedge (j \circ \rmg)) \ra Cone (i \wedge j)$ is the sum of the two composites:
   \be \begin{equation}
       \Sigma_{\rmS^1}^{-1}Cone(i \wedge j) {\overset {\Sigma_{\rmS^1}^{-1}Cone(i \wedge id)}  \ra} Cone(id_{Cone(\rmf)}  \wedge j) = Cone(\rmf) \wedge \rmB \cong Cone (\rmf \wedge id_{\rmB}) {\overset {Cone(id \wedge \rmg)} \ra} Cone(\rmf \wedge \rmg), \notag
   \end{equation} \ee
   \be \begin{equation}
       \Sigma_{\rmS^1}^{-1}Cone(i \wedge j) {\overset {\Sigma_{\rmS^1}^{-1}Cone(id \wedge j)}  \ra} \Sigma_{\rmS^1}^{-1} Cone(i \wedge id_{Cone(\rmg)}  )  = \rmA \wedge Cone(g)  \cong Cone (id _{\rmA} \wedge g) {\overset {Cone(f \wedge id)} \ra} Cone(\rmf \wedge \rmg). \notag
\end{equation} \ee
\end{lemma}
\vskip .2cm
\begin{proposition}(Multiplicative property of the pre-transfer and trace)
	\label{mult.prop}
 Assume $\rmF_i$, $i=1,2$ are simplicial presheaves, and 
 let $\rmf_i: \rmF_i \ra \rmF_i$, $i=1,2$ denote a given map. Let $\rmF = \rmF_{1+} \wedge \rmF_{2+}$
 and let $\rmf= \rmf_{1+} \wedge \rmf_{2+}$.
 Then 
 \[ tr'_{\rmF}(\rmf) = tr'_{\rmF_{1+}}(\rmf_{1+}) \wedge tr'_{\rmF_{2}}(\rmf_{2+}), \mbox { and } \]
 \[ \tau_	{\rmF}(\rmf) = \tau_{\rmF_{1+}}(\rmf_{1+}) \wedge \tau_{\rmF_{2+}}(\rmf_{2+}).   \]
 A corresponding result holds if $\rmF_2$ is a pointed simplicial presheaf with $\rmF = \rmF_{1+} \wedge \rmF_2$.
\end{proposition}
\begin{proof} A key point to observe is that the evaluation $e_{\rmF}: \rmD(\rmF) \wedge \rmF \ra \mbS_{\k}$ is given by starting with $e_{\rmF_{1+}} \wedge e_{\rmF_{2+}}: \rmD(\rmF_{1+}) \wedge \rmF_{1+} \wedge \rmD(\rmF_{2+}) \wedge \rmF_{2+} \ra \mbS_{\k} \wedge \mbS_{\k} \simeq \mbS_{\k}$ and by precomposing it
with the map $\rmD(\rmF) \wedge \rmF = \rmD(\rmF_{1+} \wedge \rmF_{2+}) \wedge \rmF_{1+} \wedge \rmF_{2+} {\overset {\tau} \ra} \rmD(\rmF_{1+}) \wedge {\rmF}_{1+} \wedge \rmD(\rmF_{2+})\wedge \rmF_{2+} $, where $\tau$ is the
obvious map that interchanges the factors. Similarly the co-evaluation map
$c: \mbS_{\k}\simeq \mbS_{\k}\wedge \mbS_{\k}{\overset {c_{\rmF_{1+}} \wedge c_{\rmF_{2+}}} \ra}\rmF_{1+}\wedge \rmD(\rmF_{1+}) \wedge \rmF_{2+} \wedge D(\rmF_{2+})$ provides the co-evaluation map for $\rmF$. The multiplicative 
property of the pre-transfer follows readily from the above two observations as well as from the definition of the pre-transfer as in Definition ~\ref{pretransfer}.
In view of the definition of the trace as in Definition ~\ref{pretransfer}, the multiplicative property of the trace follows from the
 multiplicative property of the pre-transfer. These prove the statements when $\rmF_+ = \rmF_{1+} \wedge \rmF_{2+}$. The corresponding statements when
 $\rmF_2$ is already a pointed simplicial presheaf may be proven along entirely similar lines.
\end{proof}
\vskip .2cm
\subsection{ Additivity of the Motivic Trace}
 The goal of this section is to prove the following theorem.
 \begin{theorem} (Mayer-Vietoris and Additivity for the Trace)
\label{add.trace} 
\begin{enumerate}[\rm(i)]
\item Let $\rmX$ denote a smooth quasi-projective scheme and let $i_j:\rmX_j \ra \rmX$, $j=1,2$ denote the open immersion of two Zariski open subschemes of $\rmX$, with 
$\rmX= \rmX_1 \cup \rmX_2$. Let $ \rmU \ra \rmX$ denote the open immersion of a 
Zariski open subscheme of $\rmX$, with $\rmU_i = \rmU \cap \rmX_i$. Then adopting the 
 terminology above, (that is, where $\tau_{\rmP}$ denotes the  trace associated to the pointed simplicial presheaf $\rmP$), and when $char (\k)= {\rm 0}$,
 \be \begin{equation}
 \label{MV.0}
 \tau_{\rmX/\rmU} = \tau_{\rmX_1/\rmU_1} + \tau_{\rmX_2/\rmU_2} - \tau_{(\rmX_1 \cap \rmX_2)/(\rmU_1 \cap \rmU_2)} \, in \, {\mathcal S}{\mathcal H}(k).
 \end{equation} \ee
 In case $char(\k) = \rmp > {\rm 0}$,
 \be \begin{equation}
\label{MV.p}
 \tau_{\rmX/\rmU, \mbS_{\k}[\rmp^{-1}]} = \tau_{\rmX_1/\rmU_1, \mbS_{\k}[\rmp^{-1}]} + \tau_{\rmX_2/\rmU_2, \mbS_{\k}[\rmp^{-1}]} - \tau_{(\rmX_1\cap \rmX_2)/(\rmU_1 \cap \rmU_2), \mbS_{\k}[\rmp^{-1}]}, \, in \, {\mathcal S}{\mathcal H}(k)[\rmp^{-1}].
 \end{equation} \ee
 \vskip .2cm \noindent
 Throughout the following discussion, let $<-1>$ denote the class in the Grothendieck-Witt ring associated to $-1 \in \k$
  as in \cite[p. 252]{Mo4}.
\item Let $i: \rmZ \ra \rmX$ denote a closed immersion of smooth schemes 
 with $j: \rmU \ra \rmX$ denoting the corresponding open complement. Let $\cN$ denote the normal bundle associated
 to the closed immersion $i$ and let $\Th(\cN)$ denotes its Thom-space. Let $c$ denote the codimension of $\rmZ$ in $\rmX$. Then adopting the 
 terminology above, we obtain in ${\mathcal S}{\mathcal H}(k)$ when $char (\k)= {\rm 0}$:
\be \begin{equation}
  \label{add.1.0}
\tau_{\rmX_+} = \tau_{\rmU_+} + \tau_{\rmX/\rmU},  \mbox{ and } \tau_{\rmX/\rmU} = \tau_{\Th(\cN)} = <-1>^c\tau_{\rmZ_+}.
\end{equation} \ee

In case \mbox{ $\sqrt{-1} \in k$}, it follows that 
\[\tau_{\rmX/\rmU} = \tau_{\Th(\cN)} = \tau_{\rmZ_+}.\]
  In case $char(\k) = \rmp > {\rm 0}$, we obtain in ${\mathcal S}{\mathcal H}(k)[\rmp^{-1}]$:
\be \begin{equation}
\label{add.1.p}
  \tau_{\rmX_+,  \mbS_{\k}[\rmp^{-1}]} = \tau_{\rmU_+,  \mbS_{\k}[\rmp^{-1}] } + \tau_{\rmX/\rmU, \mbS_{\k}[\rmp^{-1}]},   \tau_{\rmX/\rmU, \mbS_{\k}[\rmp^{-1}]} = \tau_{\Th(\cN), \mbS_{\k}[\rmp^{-1}]} = <-1>^c\tau_{\rmZ_+, \mbS_{\k}[\rmp^{-1}]}, 
 \end{equation} \ee
\mbox{ and assuming $\sqrt{-1} \in k$}
\[ \tau_{\rmX/\rmU, \mbS_{\k}[\rmp^{-1}]} = \tau_{\Th(\cN), \mbS_{\k}[\rmp^{-1}]} = \tau_{\rmZ_+, \mbS_{\k}[\rmp^{-1}]}.\]
\item Let $\{\rmS_{\alpha}|\alpha\}$ denote a stratification of the smooth scheme $\rmX$ into finitely many locally closed and smooth subschemes 
 $\rmS_{\alpha}$. Let $c_{\alpha}$ denote the codimension of $\rmS_{\alpha}$ in $\rmX$.  Then  we obtain in ${\mathcal S}{\mathcal H}(k)$ when $char (\k)= {\rm 0}$:
\be \begin{equation}
\label{add.2.0}
\tau_{\rmX_+} = \Sigma_{\alpha} <-1>^{c_{\alpha}}\tau_{\rmS_{\alpha+}} \mbox{ and assuming } \sqrt{-1} \in \k, 
\end{equation} \ee
\[ \tau_{\rmX_+} = \Sigma_{\alpha}\tau_{\rmS_{\alpha+}}.\]
 In case $char(\k) = \rmp > {\rm 0}$,  we obtain in ${\mathcal S}{\mathcal H}(\k)[\rmp^{-1}]$:
\be \begin{equation}
\label{add.2.p}
 \tau_{\rmX_+, \mbS_{\k}[\rmp^{-1}]} = \Sigma_{\alpha}<-1>^{c_{\alpha}} \tau_{\rmS_{\alpha+}, \mbS_{\k}[\rmp^{-1}]}, \mbox{ and again assuming }\sqrt{-1} \in \k,
 \end{equation} \ee
\[  \tau_{\rmX_+, \mbS_{\k}[\rmp^{-1}]} = \Sigma_{\alpha} \tau_{\rmS_{\alpha+}, \mbS_{\k}[\rmp^{-1}]}.\]
\end{enumerate} \qed
 \end{theorem}

\vskip .2cm \noindent
\begin{proof} We will explicitly discuss only the case in characteristic $0$, as proofs in positive characteristics will follow 
 along the same lines.
\vskip .1cm
First one recalls the stable homotopy cofiber sequence (see \cite[p. 115, Theorem 2.23]{MV})
  \be \begin{equation}
    \label{loc.0}
    \Sigma_{\T}^{\infty} \rmU_+ \ra   \Sigma_{\T}^{\infty}\rmX_+ \ra   \Sigma_{\T}^{\infty}(\rmX/\rmU) \simeq  \Sigma_{\T}^{\infty}\wedge\Th(\cN)
  \end{equation} \ee
in the stable motivic homotopy category over the base scheme. The first statement in ~\eqref{add.1.0} follows by applying Theorem ~\ref{additivity.tr.0} to the
 stable homotopy cofiber sequence in ~\eqref{loc.0}.  
 \vskip .1cm
 Next we will consider (i), namely the Mayer-Vietoris sequence. For this,  one begins with the stable cofiber sequences 
 \[\Sigma_{\T}^{\infty} (\rmU_1\cap \rmU_2)_+ \ra \Sigma_{\T}^{\infty} (\rmU_1 \sqcup \rmU_2)_+ \ra \Sigma_{\T}^{\infty} (\rmU)_+, \quad \Sigma_{\T}^{\infty} (\rmX_1\cap \rmX_2)_+ \ra \Sigma_{\T}^{\infty} (\rmX_1 \sqcup \rmX_2)_+ \ra \Sigma_{\T}^{\infty} (\rmX)_+. \] 
Then one applies Theorem ~\ref{additivity.tr.1}(i) to both of them, which will prove:
\be \begin{align}
 \label{add.2}
 \tau_{\rmU_+} =  \tau_{(\rmU_1 \cup \rmU_2)_+}  &=\tau_{\rmU_{1+}} + \tau_{\rmU_{2+}} - \tau_{(\rmU_1 \cap \rmU_2)_+} \mbox{ and } \\
 \tau_{\rmX_+} =  \tau_{(\rmX_1 \cup \rmX_2)_+}  &=\tau_{\rmX_{1+}} + \tau_{\rmX_{2+}} - \tau_{(\rmX_1 \cap \rmX_2)_+}. \notag
\end{align} \ee
On applying  the first statement in (ii)  to $\rmU_i \subseteq \rmX_i$, $i=1, 2$ and $\rmU_1 \cap \rmU_2 \subseteq \rmX_1 \cap \rmX_2$ we obtain: 
\[\tau_{\rmX_i/\rmU_i} = \tau_{\rmX_{i+}} - \tau_{\rmU_{i+}}, i=1, 2 \mbox{ and }\]
\[\tau_{(\rmX/\rmU)_+}= \tau_{(\rmX_1 \cup \rmX_2)/(\rmU_1 \cup \rmU_2)} = \tau_{(\rmX_1 \cup \rmX_2)_+} - \tau_{(\rmU_1 \cup \rmU_2)_+}.\]
The required statement in ~\eqref{MV.0} now follows on substituting from ~\eqref{add.2}.  This completes the proof of  (i).

\vskip .2cm
We proceed to establish the remaining statement in ~\eqref{add.1.0}. 
First we will consider the case where the normal bundle $\cN$ is trivial, mainly because this is an important special case to consider.
When the normal bundle is trivial, we observe that $\rmX/\rmU \simeq \Th(\cN) \simeq \T^c \wedge \rmZ_+$. Next, 
the multiplicative property of the trace as in Lemma ~\ref{mult.prop} shows that 
\be \begin{equation}
  \label{mult.prop.1}
  \tau_{\Sigma_{\T}^{\infty}(\T^c \wedge \rmZ_+)}  = (\tau_{\Sigma_{\T}^{\infty}\T})^{\wedge ^c} \wedge \tau_{\Sigma_{\T}^{\infty} \rmZ_+} 
\end{equation} \ee
as classes in $\pi_{0,0}(\mbS_{\k})$. In general, it is known that the class of $\tau_{\Sigma_{\T}^{\infty}\T}= <-1>$ in the Grothendieck-Witt group ${\rm GW}(\k)$: recall
that ${\rm GW}(\k)$ identifies with $\pi_{0,0}(\mbS_{\k})$, in view of
\cite[Theorem 6.2.2]{Mo4}. (Here it may be important to recall that $\T$ is the pointed simplicial presheaf ${\mathbb P}^1$ pointed by $\infty$.)
This implies that  $\tau_{\Sigma_{\T}^{\infty}\T} = -1  $ in $\pi_{0,0}(\mbS_{\k})$ and proves the second statement in ~\eqref{add.1.0} when $\cN$ is trivial.
\vskip .1cm
Next we assume that $ \sqrt -1 \in \k$. Then the quadratic form $<-1>$ gets identified with $<1>$ in the Grothendieck-Witt group ${\rm GW}(\k)$: see, for example, \cite[p. 44]{Sz}.
Therefore, $\tau_{\Sigma_{\T}^{\infty} \T}= <1>$, hence  $\tau_{\Sigma_{\T}^{\infty}\T^c \wedge \rmZ_+}= \tau_{\Sigma_{\T}^{\infty}\rmZ_+}$.
This completes the proof of (ii), when the normal bundle to $\rmZ$ in $\rmX$ is trivial.
\vskip .2cm
To consider the general case when the normal bundle $\cN$ is not necessarily trivial, one takes 
a finite Zariski open cover $\{\rmU_i|i=1, \cdots, n\}$ so that $\cN_{|\rmU_i}$ is trivial for each $i$. 
Then the Mayer-Vietoris property considered in (i) and ascending induction on $n$, together with the 
case where the normal bundle is trivial considered above, completes the proof in this case. (Observe that any scheme $\rmZ$ over $k$ 
of finite type is always quasi-compact, so that such a finite open cover always exists.)
 These complete the proof of all the statements in (ii). 
\vskip .2cm
Next we consider the statement in (iii).
This will follow from the  second statement in (ii) using ascending induction on the number of strata. 
However, as this induction needs to be handled carefully, we proceed to provide an outline of the
relevant argument. We will assume that the stratification of $\rmX$ defines the following increasing
 filtrations:
 \vskip .1cm
(a) $\phi=\rmX_{-1} \subseteq \rmX_0 \subseteq \cdots \subseteq \rmX_n =\rmX$, where each $\rmX_i$ is closed
and the strata $\rmX_i - \rmX_{i-1}$, $i=0, \cdots, n$ are smooth.
\vskip .1cm
(b) $\rmU_0 \subseteq \rmU_1 \subseteq \cdots \subseteq \rmU_{n-1} \subseteq \rmU_n =\rmX$, where
each $\rmU_i$ is open in $\rmX$ (and therefore smooth), with $\rmU_i - \rmU_{i-1} = \rmX_{n-i} - \rmX_{n-i-1}$, for all
$i=0, \cdots n$. Now observe that each $\rmU_k \ra \rmX$ is an 
 open immersion, while each $\rmX_k - \rmX_{k-1} \ra \rmX- \rmX_{k-1}$ is a closed immersion. Let $c_k$ denote the corresponding codimension.
\vskip .1cm
We now apply Theorem ~\ref{add.trace}(ii) with $\rmU = \rmU_{n-1}$, and $\rmZ = \rmU_n - \rmU_{n-1} = \rmX_0 - \rmX_{-1} = \rmX_0$, the closed stratum.
Since $\rmX$ is now smooth and so is $\rmZ$, the hypotheses of Theorem ~\ref{add.trace}(ii) 
are satisfied. This provides us
\be \begin{equation}
  \label{step0.iterative.add}
  \tau_{\rmX_+} = \tau_{\rmU_{n-1+}} + \tau_{\rmX/\rmU_{n-1}}\mbox{ and } \tau_{\rmX/\rmU_{n-1}} = <-1>^{c_0}\tau_{\rmX_{0+}}
\end{equation} \ee
\vskip .1cm 
Next we replace $\rmX$ by $\rmU_{n-1}$, $\rmU$ by $\rmU_{n-2}$ and $\rmZ$ by $\rmX_1- \rmX_0$. 
Since $\rmX_1 - \rmX_0$ is smooth, Theorem ~\ref{add.trace}(ii)  now provides us
\be \begin{equation}
  \label{step3.iterative.add}
 \tau_{\rmU_{n-1+}} = \tau_{\rmU_{n-2+}} + <-1>^{c_1} \tau_{(\rmX_1 - \rmX_0)_+}.
\end{equation} \ee
\vskip .1cm
Substituting these in ~\eqref{step0.iterative.add}, we obtain
\[ \tau_{\rmX_+} = \tau_{\rmU_{n-2+}} + <-1>^{c_1} \tau_{(\rmX_1- \rmX_0)_+} + <-1>^{c_0}\tau_{\rmX_{0+}}. \]
Clearly this may be continued inductively to deduce statement (iii) in Theorem ~\ref{add.trace}
 from Theorem ~\ref{add.trace}(ii). 
 \end{proof}
 \section{\bf Proofs of the main Theorems}
 
 We begin by discussing the following Proposition, which seems to be rather well-known. (See for example, \cite[Proposition 4.10]{Th86} or
 \cite[(3.6)]{BP}.)
\begin{proposition}
 \label{struct.T.actions}
Let $\rmT$ denote a split torus acting on a separated scheme $\rmX$ all defined over the given perfect base field $k$. 
\vskip .1cm
Then the following hold.
\vskip .1cm
$\rmX$ admits a  decomposition into a disjoint union of finitely many locally closed, $\rmT$-stable subschemes $\rmX_j$ so that
 \be \begin{equation}
    \label{torus.decomp}
    \rmX_j \cong (\rmT/\Gamma _j) \times \rmY_j.
 \end{equation} \ee
 \vskip .2cm \noindent
 Here each $\Gamma_j$ is a sub-group-scheme of $\rmT$, each $\rmY_j$ is a scheme of finite type over $k$ which is also regular and on which $\rmT$ acts trivially with the isomorphism in ~\eqref{torus.decomp} being $\rmT$-equivariant. 
\end{proposition}
\begin{proof}
One may derive this from the generic torus slice theorem proved in \cite[Proposition 4.10]{Th86}, which says that if a split torus acts on a reduced separated  scheme of finite type over a perfect field, then the following are satisfied:
 \begin{enumerate}
   \item there is an open subscheme $\rmU$ which is regular and stable under the $\rmT$-action
   \item a geometric quotient $\rmU/\rmT$ exists, which is a regular scheme of finite type over $k$
   \item $\rmU $ is isomorphic as a $\rmT$-scheme to $\rmT/\Gamma \times \rmU/\rmT$ where $\Gamma$ is a diagonalizable subgroup scheme of $\rmT$ and $\rmT$ acts trivially on $\rmU/\rmT$. 
 \end{enumerate}
(See also \cite[(3.6)]{BP} for a similar decomposition.) 
\end{proof}
Next we consider the following theorem.
 \begin{theorem}
  \label{torus.act} 
  Under the assumption that the base field $k$ is of characteristic $0$,  the following hold, where $\tau_{\rmX_+}$ denotes the trace associated to the pointed scheme $\rmX_+$:
  \begin{enumerate}[\rm(i)]
      \item $\tau_{{\mathbb G}_{m+}}=1 -<-1>$ in ${\rm GW}(\k)$, and if $\rmT$ is a split torus of rank $n$, $\tau_{\rmT_+}= (1-<-1>)^n$  in ${\rm GW}(\k)$.
      Therefore, it follows that when $\k$ contains a $\sqrt -1$, $\tau_{{\mathbb G}_{m+}} =0 $ and $\tau_{\rmT_+} = 0$ in ${\rm GW}(\k)$.
      \item Let $\rmT$ denote a split torus acting on a smooth scheme $\rmX$. Then $\rmX^{\rmT}$ is also 
            smooth, and $\tau_{\rmX_+}- \tau_{\rmX^{\rmT}_+}$ belongs to the ideal generated by $(1- <-1>)$ in ${\rm GW}(\k)$. In particular, when
            $\k$ contains a $\sqrt -1$, $\tau_{\rmX_+}= \tau_{\rmX^{\rmT}_+}$ in ${\rm GW}(\k)$.
  \end{enumerate}      
 \vskip .1cm \noindent
    If the base field is of positive characteristic $\rmp$, the corresponding assertions hold
      with the the trace of a pointed smooth scheme $\rmY_+$ replaced by $\tau_{\rmY_+, \mbS_{\k}[\rmp^{-1}]}$ and the
      Grothendieck-Witt ring  replaced by the Grothendieck-Witt ring with the prime $\rmp$ inverted.
  \vskip .2cm \noindent  
 \end{theorem}
 \begin{proof} We will only consider the proofs when the base field is of characteristic $0$, since the proofs in the positive characteristic case are entirely similar. However, it is important to point out that in positive characteristics $\rmp$, it is important to invert $\rmp$: 
 for otherwise, one no longer has a theory of Spanier-Whitehead duality.
 Next observe from Definition ~\ref{pretransfer}, that the trace $\tau_{\rmX_+}$ associated to any smooth scheme $\rmX$ is a map
 $\mbS_{\k} \ra \mbS_{\k}$: as such, we will identify $\tau_{\rmX_+}$ with the corresponding class $\tau_{\rmX_+}^*(1)$ in the Grothendieck Witt-ring of the base field.
 
 \vskip .1cm 	
 Next we consider (i). We observe that the scheme ${\mathbb A}^1$ is the disjoint union of 
 the closed point $\{0\}$ and ${\mathbb G}_m$. If $i_1: \{0\} \ra {\mathbb A}^1$ and $j_1: {\mathbb G}_m \ra {\mathbb A}^1$ are the corresponding immersions, Theorems ~\ref{add.trace}(ii) and (iii) 
  show that
  \be\begin{equation}
 	  \label{torus.1.0}
 	   \tau_{{\mathbb A}^1_+} =   \tau_{{\mathbb G}_{m+}} + \tau_{{\mathbb A}^1/{\mathbb G}_{m}} =  \tau_{{\mathbb G}_{m+}} + \tau_{\T} =  \tau_{{\mathbb G}_{m+}} + <-1>. 
 	\end{equation} \ee
Therefore, it follows that 
 	\be\begin{equation}
 	  \label{torus.1.1}
 	    \tau_{{\mathbb G}_{m+}} = \tau_{{\mathbb A}^1_+} - <-1>  = 1-<-1>. 
 	\end{equation} \ee
 \vskip .1cm \noindent
 where $\tau_{{\mathbb A}^1_+} = \tau_{\{0\}_+} =1$ by ${\mathbb A}^1$-contractibility. One may readily see this from the
 definition of the pre-transfer as in Definition ~\ref{pretransfer}, which shows that both the
 pre-transfer $\tr_{\rmC_+}' = \tr_{\rmC_+}'(id)$ and hence the corresponding trace, $\tau_{\rmC_+} = \pi \circ tr_{\rmC_+}'$ depends on $\rmC_+$ only up to its class in the motivic stable homotopy category.
 Since $\rmT$ is a split torus, we may assume $\rmT = {\mathbb G}_{\it m} ^n$ for some 
 positive integer $n$. Then the multiplicative property of the trace and pre-transfer (see Proposition ~\ref{mult.prop}) prove that
  $\tau_{\rmT_+}= (1-<-1>)^n$. In particular, when $\k$ contains a $\sqrt -1$, it follows that $\tau_{{\mathbb G}_{m+}} =0$ and
  $\tau_{\T+}=0$ in ${\rm GW}(\k)$. These complete the proof of statement (i).
 \vskip .2cm
 Therefore, we proceed to prove the statement in (ii). First, we invoke Proposition ~\ref{struct.T.actions}
to conclude that $\rmX^{\rmT}$ is the disjoint union of the subschemes $\rmX_j$ for which $\Gamma _j = \rmT$. 
 \vskip .1cm
 Let $i_j: \rmX_j \cong (\rmT/\Gamma _j) \times \rmY_j \ra \rmX$ denote the locally closed immersion.
 Next observe that the additivity of the trace proven in Theorem ~\ref{add.trace},  and the multiplicativity of the pre-transfer and trace proven in  
 Proposition ~\ref{mult.prop} along with the decomposition in  ~\eqref{torus.decomp} show that
 \be \begin{align}
    \label{tr.T}
    \tau_{\rmX_+} = \Sigma _j \tau_{\rmX_{j+}} &= \Sigma _j  (\tau_{\rmT/\Gamma_j+}) \wedge \tau_{\rmY_{j+}}.
 \end{align} \ee
 \vskip .2cm
 Now statement (i)  in the theorem shows that the term $\tau_{\rmT/\Gamma_j+} = (1-<-1>)^{n_j}$, if $\rmT/\Gamma_j$ is a split torus of rank $n_j$.
 Since $\rmX^{\rmT}$ is the disjoint union of the subschemes $\rmX_j = \rmT/\Gamma_j \times \rmY_j$ with $\Gamma_{\rm j}= \rmT$, the additivity of the trace proven in 
 Theorem ~\ref{add.trace} and applied to $\rmX^{\rmT}$ proves the sum of such terms on the right-hand-side of ~\eqref{tr.T} is 
 $\tau_{\rmX^{\rmT}_+}$. Therefore, it follows that $\tau_{\rmX_+} - \tau_{\rmX^{\rmT}_+}$ belongs to the ideal in ${\rm GW}(\k)$ generated by $1-<-1>$.
 In particular, when $\k$ contains a $\sqrt -1$, it follows that $\tau_{\rmX+} = \tau_{\rmX^{\rmT}_+}$. These complete the proof of the statements in (ii).
\end{proof}
 
 \vskip .2cm \noindent
 {\bf Proof of Theorem ~\ref{main.thm}.} We point out that it is important to assume the base field $k$ is perfect in the following arguments: this will ensure
 that all the schemes considered here are defined over the same base field.
 First we will show that we can reduce to the case $\rmG$ is {\it connected}. Let $\rmG^o$ denote the connected
 component of $\rmG$ containing the identity element and let $\rmT$ denote a split maximal torus in $\rmG$. Then,  one first obtains the
 isomorphisms $\rmG/\rmN_{\rmG}(\rmT) \cong \{\g\rmT \g^{\rm-1}|\g \eps \rmG\}$ and $\rmG^o/\rmN_{\rmG^o}(\rmT) \cong \{\g_o\rmT \g_o^{\rm-1}|\g_o \eps \rmG^o\}$.
 Next observe that $\g\rmT \g^{\rm -1}$, being a maximal torus and hence a connected subgroup of $\rmG$, is in fact a maximal torus in $\rmG^o$ for 
 each $\g \eps \rmG$. These show that 
 \[\rmG/\rmN_{\rmG}(\rmT) = \{\g\rmT \g^{\rm -1}|\g \eps \rmG\} \cong \{\g_o\rmT \g_o^{\rm -1}|\g_o \eps \rmG^o\} = \rmG^o/\rmN_{\rmG^o}(\rmT).\]
 Therefore, we may assume the group $\rmG$ is connected.
 \vskip .1cm
Moreover, we may take the quotient by the unipotent radical $\rmR_u(\rmG)$, which is a normal subgroup
  (and is isomorphic to an affine space), with the quotient $\rmG_{red} = \rmG/\rmR_u(\rmG)$ reductive. 
  Now $\rmG/\rmN_{\rmG}(\rmT) \cong \rmG_{red}/\rmN_{\rmG_{red}}(\rmT)$ (since the intersection of a maximal torus in $\rmG$ with the unipotent radical $\rmR_u(\rmG)$ is trivial), so that we may assume 
  $\rmG$ is a connected split reductive group.
 \vskip .1cm
 Then we observe that since $\rmG/\NT$ is the 
variety of all split maximal tori in $\rmG$,  $\rmT$ has an action on $\rmG/\NT$ (induced by the left translation action of $\rmT$ on $\rmG$) 
so that there is exactly a single fixed point, namely the coset $e\NT$, that is, $(\rmG/\NT)^{\rmT} = \{e\NT\} = \{Spec \, \k \}$.
(To prove this assertion, one may reduce to the case where the base field is algebraically closed, since the formation
of fixed point schemes respects change of base fields as shown in \cite[p. 33, Remark (3)]{Fog}. See also \cite[Lemma 3.5]{BP}.
In fact, one may see this directly as follows. Making use of the identification of $\rmG/\rmN_{\rmG}(\rmT)$ with $\{\g\rmT \g^{\rm-1}|\g \eps \rmG\}$, one sees that if $g_0\rmT {\it g}_0^{-1}$ is fixed by the conjugation action of $\rmT$, 
then $g_0^{-1}\rmT {\it g}_0 \subseteq \rmN_{\rmG}(\rmT)^o = \rmT$, so that $g_0 \eps \rmN_{\rmG}(\rmT)$. Thus the coset $g_0\rmN_{\rmG}(\rmT) = e\rmN_{\rmG}(\rmT)$.)
\vskip .2cm
Next we will first consider the case where the base field is of characteristic $0$. Therefore, by Theorem ~\ref{torus.act}(ii),
\[\tau_{\rmG/\NT_+} = \tau_{(\rmG/\NT)^{\rmT}_+} =\tau_{{{\rm Spec}\, k}_+} = id_{\mbS_{\k}},\]
which is the identity map of
the motivic sphere spectrum. Therefore, 
\[\chi_{mot}(\rmG/\NT)=\tau_{\rmG/\NT_+}^*(1) =1.\]
The motivic stable homotopy group $\pi_{0, 0}(\mbS_{\k})$
identifies with the Grothendieck-Witt ring by \cite{Mo4}. This completes the proof of the statement on $\tau_{\rmG/\NT_+}$ in
Theorem ~\ref{main.thm} in this case. In case the base field is of positive characteristic $\rmp$,
one observes that $\Sigma_{\T}^{\infty}\rmG/\NT_+$ will be dualizable only in $\Spt(\k_{mot})[{\rm p}^{\rm -1}]$.
But once the prime $\rmp$ is inverted the same arguments as before carry over proving the corresponding statement.
These complete the proof of the theorem. \qed
\vskip .2cm \noindent
{\bf Proof of Corollary ~\ref{cor}}.
\vskip .1cm
Observe, first that if $\bar k$ is the algebraic closure of the given field,
then it contains a $\sqrt -1$, and therefore the  conclusions of the theorem hold in this case. In positive characteristic $\rmp$, we proceed to show that this already implies that $\chi_{mot}(\rmG/\NT)$ is a {\it unit}
	in the group ${\rm GW}({\it k})[\rmp^{-1}]$, 
	without the assumption on the existence of a square root of $-1$ in ${\it k}$. For this, one may first observe the commutative diagram, where
	$\bar {\it k}$ is an algebraic closure of ${\it k}$:
	\be \begin{equation}
	     \label{GW.pos.char}
	    \xymatrix{ { {\rm GW}( \bar {\it k})[\rmp^{-1}]} \ar@<1ex>[r]^(.6){rk}_(.6){\cong} & {{\mathbb Z}[\rmp^{-1}]}\\
	               { {\rm GW}( {\it k})[\rmp^{-1}]} \ar@<1ex>[r]^(.6){rk} \ar@<1ex>[u] \ar@<1ex>[r] & {{\mathbb Z}[\rmp^{-1}]} \ar@<1ex>[u]^{id}.}
	    \end{equation}
Here the left vertical map is induced by the change of base fields from ${\it k}$ to $ \bar {\it k}$, and $rk$ denotes the {\it rank} map. Since the motivic Euler-characteristic
of $\rmG/\NT$ over $Spec \, {\it k}$ maps to the motivic Euler-characteristic of the corresponding $\rmG/\NT$ over $Spec \, \bar {\it k}$,
it follows that the rank of $\chi_{mot}(\rmG/\NT)$ over $Spec \, {\it k}$ is in fact $1$. By \cite[Lemma 2.9(2)]{An}, this shows that the
 $\chi_{mot}(\rmG/\NT)$ over $Spec \, {\it k}$ is in fact a unit in ${\rm GW}( {\it k})[\rmp^{-1}]$, that is , 
 when ${\it k}$ has positive characteristic. (For the convenience of the reader, we 
 will summarize a few key facts discussed in \cite[Proof of Lemma 2.9(2)]{An}. It is observed there that when the base field $k$ is {\it not} formally real,
 then ${\rm I}(k) = kernel({\rm GW}(k) {\overset {rk} \ra }{\mathbb Z})$ is the nil radical of ${\rm GW}(k)$: see \cite[Theorem V.8.9, Lemma V.7.7 and Theorem V. 7.8]{Bae}. Therefore,
 if $char(k)= \rmp> {\rm 0}$, and the rank of $\chi_{mot}(\rmG/\NT)$ is $1$ in ${{\mathbb Z}[\rmp^{-1}]}$, then $\chi_{mot}(\rmG/\NT)$ is $1+q$ for some nilpotent element $q$ in ${\rm I}(k)[\rmp^{-1}]$ and 
 the conclusion follows.)
 \vskip .1cm
 An alternative shorter proof is the following: observe that $ \chi_{mot}(\rmG/\NT) - \chi_{mot}((\rmG/\NT)^{\rmT}) = \chi_{mot}(\rmG/\NT) -1 $ belongs
 to the ideal generated by $1-<-1>$. $1 - <-1>$ clearly belongs to ${\rm I}(k) = kernel({\rm GW}(k) {\overset {rk} \ra }{\mathbb Z})$. When 
 $\k$ is not formally real, the above ideal is nilpotent as observed above, and therefore, $\chi_{mot}(\rmG/\NT)$ is a unit when $\k$ is not formally real.
\vskip .1cm
In characteristic $0$, the commutative
diagram 
\be \begin{equation}
	     \label{GW.pos.char}
	    \xymatrix{ { {\rm GW}( \bar {\it k})} \ar@<1ex>[r]^(.6){rk}_(.6){\cong} & {{\mathbb Z}}\\
	               { {\rm GW}( {\it k})} \ar@<1ex>[r]^(.6){rk} \ar@<1ex>[u] \ar@<1ex>[r] & {{\mathbb Z}} \ar@<1ex>[u]^{id}.}
	    \end{equation}
shows that once again the rank of $\chi_{mot}(\rmG/\NT)$ is $1$. Therefore,
 to show that the class\\ $\chi_{mot}(\rmG/\NT)$ is a unit in ${\rm GW}( {\it k})$, it suffices to show its signature is 
 a unit: this is proven in \cite[Theorem 5.1(1)]{An}. (Again, for the convenience of the reader,
 we summarize some details from the proof of \cite[Theorem 5.1(1)]{An}. When the field $k$ is not formally real, the discussion in
 the last paragraph applies, so that by \cite[Lemma 2.12]{An} one reduces to considering only the case when $k$ is a real closed field. In this case, one lets ${\mathbb R}^{alg}$ denote 
 the real closure of ${\mathbb Q}$ in ${\mathbb R}$. Then, one knows the given real closed field $k$ contains a copy of ${\mathbb R}^{alg}$ and that there exists
 a reductive group scheme $\widetilde \rmG$ over  ${\rm Spec \, } {\mathbb R}^{alg}$ so that $\rmG = \widetilde \rmG \times_{{\rm Spec \, } {\mathbb R}^{alg}} {\rm Spec \, }k$. Let
 $\rm G_{\mathbb R} = \widetilde \rmG \times _{{\rm Spec \, } {\mathbb R}^{alg}} {\rm Spec \, }{\mathbb R}$.  Then 
 one also  observes 
 that the Grothendieck-Witt groups of the three fields $k$, ${\mathbb R}^{alg}$ and ${\mathbb R}$ are isomorphic, and the
 motivic Euler-characteristics $\chi_{mot}(\rmG/\NT)$,  $\chi_{mot}({\widetilde \rmG}/ {\widetilde {\NT }})$ and
 $\chi_{mot}({\rmG}_{{\rm Spec \, }{\mathbb R}} / { \rmN(T) }_{{\rm Spec \, }{\mathbb R}})$
 over the above three fields identify under the above isomorphisms, so that one may assume 
 the base field $k$ is ${\mathbb R}$. Then it is shown in \cite[Proof of Theorem 5.1(1)]{An} that, in this case,  knowing the rank and signature of the
 motivic Euler characteristic $\chi_{mot}(\rmG/\NT)$ are $1$ suffices to 
 prove it is a unit in the Grothendieck-Witt group.)
 These complete the proof of the corollary. \qed

 \vskip .3cm
 
 
 \end{document}